
\magnification=\magstep1
\input amstex
\documentstyle{amsppt}
\leftheadtext{E. Makai, Jr.* , H. Martini, T. \'Odor}
\rightheadtext{On a theorem of Ryabogin and Yaskin}
\topmatter
\title On a theorem of D. Ryabogin and 
V. Yaskin about detecting symmetry\endtitle
\author E. Makai, Jr.*, H. Martini, T. \'Odor\endauthor
\address Alfr\'ed R\'enyi Mathematical Institute, 
Hungarian Academy of Sciences,
\newline
H-1364 Budapest, Pf. 127, Hungary
{\rm{http://www.renyi.mta.hu/\~{}makai}}
\vskip.1cm
Fakult\"at f\"ur Mathematik, Technische Universit\"at Chemnitz
\newline
D-09107 Chemnitz, Germany
\vskip.1cm
Bolyai Institute, University of Szeged,
\newline
H-6720 Szeged, Aradi V\'ertan\'uk tere 1, Hungary
\endaddress
\email makai.endre\@renyi.mta.hu,
\newline
martini\@mathematik.tu-chemnitz.de
\newline
odor\@math.u-szeged.hu
\endemail
\thanks *Research (partially) supported by Hungarian National Foundation for 
Scientific Research, grant nos. K75016, K81146.\endthanks
\subjclass {\it 2010 Mathematics Subject Classification.} Primary:
44A99. Secondary: 52A38
\endsubjclass
\abstract We give a simple deduction 
of a recent theorem of D. Ryabogin and
V. Yaskin, about detecting symmetry of star bodies in ${\Bbb{R}}^n$
with $C^1$ radial functions  --- via their conical section functions --- 
from an older theorem of us.\endabstract
\endtopmatter\document


\head \S 1 Notions and notations \endhead

We will work in {\it{Euclidean space}}\,\, ${\Bbb{R}}^n$, where 
$n \ge 2$. Its {\it{unit sphere}} will be written as $S^{n-1}$.
We say that $K \subset {\Bbb{R}}^n$ is a {\it{star body}} if it is of the
form $K=\{ \lambda u \mid u \in S^{d-1},\,\, 0 \le \lambda \le 
\varrho _K(u) \} $,
where $\varrho _K: S^{n-1} \to (0, \infty )$ is a continuous function, which is
called the
{\it{radial function}} of the star body $K$. 
A {\it{convex body in}} ${\Bbb{R}}^n$
is a compact convex set with interior points.
If $K \subset {\Bbb{R}}^n$ is a convex body 
containing $0$ in its interior, then it is a star
body. Moreover, its radial function $\varrho _K$
is Lipschitz (we consider $S^{n-1}$ with
its geodesic metric, and Lipschitz is meant with respect to it),
cf. [MM\'O1], first paragraph of \S 3. 
For $f:S^{n-1} \to {\Bbb{R}}$ a Lipschitz function, we denote by $L(f)$ its
Lipschitz constant (with respect to the geodesic metric on $S^{n-1}$).
Observe that if the radial function
$\varrho _K$ of a star body $K$ is $C^1$ (for
which we will shortly say that {\it{the star
body is}} $C^1$), then it is Lipschitz (actually, this
implication holds even for any function $S^{n-1} \to {\Bbb{R}}$). 

For $\xi \in S^{n-1}$ we write $\xi ^{\bot }$ for the linear $(n-1)$-subspace
of ${\Bbb{R}}^n$
orthogonal to $\xi $. We will use also spherical polar coordinates, with
north pole some $\xi \in S^{n-1}$. That is, we write each $x \in S^{n-1}$ as
$$
x= \xi \sin \psi + \eta \cos \psi ,{\text{ \,\,where\,\, }} 
\eta \in S^{n-1} \cap
\xi ^{\bot } {\text{ \,\,and\,\, }} -\pi /2 \le \psi \le \pi /2 .
$$ 
We call $\psi $ the {\it{geographic latitude}} (which will be more convenient
for us than the customarily used $\varphi = \pi /2 - \psi $), and will write
$$
x= (\eta , \psi ) .
$$

A function $f:S^{n-1} \to {\Bbb{R}}$ is {\it{even}} if
$f(x)=f(-x)$, for all $x \in S^{n-1}$.

D. Ryabogin and V. Yaskin 
[RY], p. 509, denoted, for \,$\xi \in S^{n-1}$ and $z \in (-1,1)$,
by $C(\xi ,z)$ the cone
$\{ 0 \} \cup \{ x \in {\Bbb{R}}^n \setminus \{ 0 \} \mid \cos (\angle \xi 0 x) 
=z \} $. 
Then, for $K \subset {\Bbb{R}}^n$ a star body, [RY], pp. 509-510,
defined the {\it{conical section function}} $C_{K, \xi } (z)$
of $K$ as
$$
C_{K,\xi } (z):= {\text{vol}}_{n-1} \left( K \cap C(\xi ,z) \right) ,
$$
where ${\text{vol}}_{n-1}$ means $(n-1)$-{\it{volume}}.


\head \S 2 Some results of D. Ryabogin-V. Yaskin and 
E. Makai, Jr.-H. Martini-T. \'Odor 
\endhead

[RY] proved the following geometrical theorem, by a relatively short proof, 
but using advanced methods, namely, Fourier transform techniques. (The
converse implication in Theorem A is obvious.)

\proclaim{Theorem A} ({\rm{[RY]}}, Theorem 1.1)
Let $K \subset {\Bbb{R}}^n$ be a \,$C^1$ star body. 
Assume that, for all \,$\xi \in S^{n-1}$, the
function $C_{K, \xi } (z)$ has a critical point at $z=0$. Then the body $K$ is
$0$-symmetric. 
\endproclaim

Here, a {\it{critical point of a function}} is a point such that
the derivative of the function
at this point exists, and equals $0$. 

[RY], in the remarks in the second paragraph after their Theorem 1.2,
mentioned that, by the methods of [MM\'O1], Theorem A can be extended to 
any convex body containing $0$ in its interior. This also follows from our
Theorem A$'$ below.


Theorem A follows from the following analytical theorem.

\proclaim{Theorem A$'$}
Let $f:S^{n-1} \to {\Bbb{R}}$ be a Lipschitz function.  
Assume that, for almost all \,$\xi \in S^{n-1}$,
 we have that the integral of $f$ 
on the set $S^{n-1} \cap C(\xi ,z)$, as a function of $z$, 
has a critical point at $z=0$. Then $f$ is an even function.
\endproclaim

We obtain Theorem A, by applying Theorem A$'$ to the $C^1$, hence
Lipschitz function $f:=\varrho _K ^{n-1} /(n-1)$. 
Since the radial
function $\varrho _K$ of a convex body $K \subset {\Bbb R}^d$, 
containing $0$ in its interior, is Lipschitz (cf. \S 1), 
the above mentioned extension of Theorem
A to convex bodies, containing $0$ in their interiors,
follows from Theorem A$'$ similarly.


\definition{Remark}
To justify the hypotheses of Theorem
A$'$, we recall from [MM\'O1], Lem\-ma 3.5 and
its proof, and Lemma 3.6,
the following. For $f:S^{n-1} \to {\Bbb{R}}$ being a
Lipschitz function, 
for almost all \,$\xi \in S^{n-1}$ we have that, for almost
all \,$x \in S^{n-1} \cap \, \xi ^{\bot }$, the function $f$ is
differentiable. 
Further, 
for almost
all \,$\xi \in S^{n-1}$ we have that, for $z=0$,
$$
\frac{d}{dz} \int _{S^{n-1} \cap (\xi ^{\bot} + z \xi) } f(x)dx
$$
exists. Moreover, it equals
$$
\int _{S^{n-1} \cap \xi ^{\bot} } \frac{\partial f}{\partial \psi } (x) dx,
$$ 
where
$\psi = \psi _{\xi }$ is the geographic latitude, with the north pole at $\xi $ 
(hence the partial derivative $\partial f / \partial \psi $
is taken along a meridian, in the direction toward the north pole $\xi $), 
and where also the second integral exists, for almost
all \,$\xi \in S^{n-1}$. (These readily 
imply that, also in Theorem A$'$, the converse implication holds.) 
\enddefinition
  

\newpage

Now we cite a theorem from [MM\'O1].

\proclaim{Theorem B} ({\rm{[MM\'O1]}}, Lemma 3.6, Theorem 3.8)
Let $f:S^{n-1} \to {\Bbb{R}}$ be a Lipschitz function.  
Assume that, for almost all \,$\xi \in S^{n-1}$, we have 
$$
\int _{S^{n-1} \cap \xi ^{\bot} }\frac{\partial f}{\partial \psi } (x) dx =0.
$$
Then $f$ is an even function.
\endproclaim

To justify the hypotheses of Theorem B, recall the Remark above. 
(The above remark readily implies that, also in Theorem B,
the converse implication holds.) 


We have used Theorem B in [MM\'O1] to prove another geometrical theorem.
This theorem was proved for $n=2$ by [H], Theorem 1; for $n \ge 3$ it was
first proved by [MM\'O1], Corollary 3.4, Lemma 3.5, Theorem 3.8,
by using spherical harmonics, and the Funk-Hecke formula.
It was reproved, for $n \ge 3$, by a relatively short proof,
however, using advanced methods, namely, Fourier transform techniques,
by [RY], Theorem 1.2, for the $C^1$ case.
This geometrical theorem states the following.
Let $K \subset {\Bbb{R}}^n$ be a star body with Lipschitz radial function.
Then, for almost all \,$\xi \in S^{n-1}$, the function 
$z \mapsto {\text{\rm{vol}}}_{n-1}\left( K \cap (\xi ^{\bot } + z \xi )
\right) $ (${\text{\rm{vol}}}_{n-1}$ meant here as $(n-1)$-dimensional
Lebesgue measure)
is differentiable at \,$0$.
Let, for almost all \,$\xi \in S^{n-1}$, this
function 
have a critical point at $z=0$. Then the body $K$ is
$0$-symmetric. (The last but two sentence readily 
implies that, also in this geometrical theorem, 
the converse implication holds.)


An infinitesimal variant of the last mentioned geometrical theorem, for the
case of a convex body (infinitesimally) close to the unit ball,
not for the $(n-1)$-volumes
of the intersections $K \cap (\xi ^{\bot } + z \xi )$, 
but for the ($(n-2)$-dimensional) surface area, and also for the lower (but
positive)
dimensional quermassintegrals (cf. [BF], \S 32, [Sch], {\S}{\S}
4.1, 4.2) of these
intersections, has been proved, 
in the sufficiently regular case, in [MM], Theorem.
Details cf. there. 


In [MM\'O2] we have proved an (almost) generalization of Theorem B, when 
$\partial f / \partial \psi $ in the hypothesis of Theorem B 
was replaced by $(\partial / \partial \psi )^m f$, for $m \ge 2$ an integer.
Details cf. there.


In what follows, we show that our Theorem B implies Theorem A$'$ (and thus also
Theorem A).


\head \S 3 Proof of the implication Theorem B $\Longrightarrow $ Theorem
A$'$.
\endhead

\demo{Proof}
We have, writing $\sin \psi :=z$ (where $ - \pi /2 < \psi < \pi /2$),
$$
\int _{S^{n-1} \cap C(\xi , z)} f(x)dx = 
\int _{S^{n-1} \cap C(\xi , z)} f (\eta , \psi )  d(\eta, \psi ) .
\tag 1
$$
By the Lipschitz property of $f$ we have, for $(\eta , \psi ) \in S^{n-1}$,
$$
|f (\eta , \psi ) - f (\eta, 0 )| \le L(f) \cdot | \psi | .
\tag 2
$$
By \thetag{2} we have, for $\psi \ne 0$, that

\newpage

$$
\cases
|f (\eta , \psi ) - f (\eta, 0 )|/ | \sin \psi |
=[|f (\eta , \psi ) - f (\eta, 0 )|/ | \psi |] \cdot [\psi / \sin \psi ] \\
< L(f) \cdot (\pi /2) / 1 =:c . 
\endcases
\tag 3
$$
By 
[MM\'O1], Lemma 3.5 and
its proof, and Lemma 3.6 (cf. also our Remark),
for almost all \,$\xi \in S^{n-1}$, we have that for almost all
$x \in S^{n-1} \cap \xi ^{\bot }$ the function $f$ is differentiable, and thus,
in particular, $(\partial f / \partial \psi )(x)$ 
(taken along a meridian, in the direction toward $\xi $)
exists. Moreover, for these 
(almost all) $\xi $'s, and for $z \in (-1,1)
\setminus \{ 0 \} $ and $z \to 0$ we have, 
using in the first equality \thetag{1}, 
$$
\cases
\left[ \int _{S^{n-1} \cap C(\xi , z)} f(x)dx - 
\int _{S^{n-1} \cap C(\xi , 0)} f(x) dx \right] \big/ z = \\
\left[ \int _{S^{n-1} \cap C(\xi , z)} f(\eta , \psi )d(\eta , \psi ) - 
\int _{S^{n-1} \cap C(\xi , 0)} f(\eta , 0)d(\eta , 0 ) \right] \big/ z = \\
\left[ \int _{S^{n-1} \cap \xi ^{\bot }} f(\eta , \psi ) d \eta \cdot \cos
^{n-2} \psi -
\int _{S^{n-1} \cap \xi ^{\bot }} f(\eta , 0 ) d \eta \right] \big/ \sin \psi =
\\
\int _{S^{n-1}\cap \xi ^{\bot }} \left[ [  f(\eta , \psi ) d \eta -
f(\eta , 0 ) ] \big/ \psi ] \cdot [ \psi / \sin \psi ] \right] d \eta
+ \\
\int _{S^{n-1} \cap \xi ^{\bot }} f(\eta , \psi ) d \eta \cdot
( \cos ^{n-2} \psi -1) \big/ \sin \psi
\to \\
\int _{S^{n-1} \cap \xi ^{\bot }} (\partial f / \partial \psi ) 
(\eta ) \cdot 1 \cdot d \eta +0,
\endcases
\tag 4
$$
by $ \psi / \sin \psi  \to 1, $ and
$|\int _{S^{n-1} \cap \xi ^{\bot }} f(\eta , \psi ) d \eta | \cdot
(1 - \cos ^{n-2} \psi ) \big/ |\sin \psi | \le 
\cos ^{n-2}( \psi ) \times $
${\text{\rm{vol}}}_{n-2}(S^{n-2}) \cdot \max \{ |f(x)|
\mid x \in S^{n-1} \} \cdot 
(1 - \cos ^{n-2} \psi ) \big/ |\sin \psi | =O(| \psi |) \to 0, $
for $ \psi \to 0$ (${\text{vol}}_{n-2}$ meaning $(n-2)$-volume). 
Still we used for the convergence of the summand in the fourth line  
of \thetag{4} 
Lebesgue's dominated convergence theorem with
integrable majorant $c$, cf. \thetag{3},
for each $\xi \in S^{n-1}$ for which for almost all
$x \in S^{n-1} \cap \xi ^{\bot }$ the function $f$ is differentiable, thus for
almost all \,$\xi \in S^{n-1}$.

By the hypothesis of Theorem A$'$, the last expression in \thetag{4} vanishes
for almost all \,$\xi \in S^{n-1}$,
thus the hypothesis of Theorem B is satisfied. Hence also the conclusion of
Theorem B is satisfied, i.e., $f$ is even, which is the conclusion of
Theorem A$'$ as well.
$ \blacksquare $

\enddemo


\enddocument